\documentclass[12pt,leqno]{article}
\usepackage{amsmath,amssymb,amscd,latexsym}
\def\bbbone{{\mathchoice {\rm 1\mskip-4mu l} {\rm 1\mskip-4mu l}
{\rm 1\mskip-4.5mu l} {\rm 1\mskip-5mu l}}}
\def \={{\ =\ }}
\newcommand{\1}{{\bbbone}}
\newcommand{\C}{{\mathbb{C}}}
\newcommand{\N}{{\mathbb{N}}}
\newcommand{\Q}{{\mathbb{Q}}}
\newcommand{\R}{{\mathbb{R}}}
\newcommand{\Z}{\mathbb{Z}}
\newcommand{\ddet}{\mathrm{det}}
\newcommand{\id}{\mathrm{id}}
\renewcommand{\mod}{\;\mathrm{mod}\;}
\newcommand{\sgn}{\mathrm{sgn}\,}
\newcommand{\Fix}{\mathrm{Fix}}
\newcommand{\GL}{\mathrm{GL}\,}
\newcommand{\Hom}{\mathrm{Hom}}
\newcommand{\Imm}{\mathrm{Im}\,}
\newcommand{\Ker}{\mathrm{Ker}\,}
\newcommand{\Lie}{\mathrm{Lie}\,}
\newcommand{\PSL}{\mathrm{PSL}}
\newcommand{\Tr}{\mathrm{Tr}}
\newcommand{\Fh}{{\mathcal F}}
\newcommand{\Hh}{{\mathcal H}}
\newcommand{\Rh}{{\mathcal R}}
\newcommand{\ea}{{\mathfrak{a}}}
\newcommand{\eg}{{\mathfrak{g}}}

\newcommand{\ep}{\mathfrak{p}}
\newcommand{\eZ}{\mathfrak{c}}
\newcommand{\oH}{\bar{H}}
\newcommand{\obeta}{\bar{\beta}}
\newcommand{\okappa}{\bar{\kappa}}
\newcommand{\ohne}{\setminus}
\newcommand{\silo}{\stackrel{\sim}{\longrightarrow}}
\newcommand{\tei}{\, | \,}
\newcommand{\hullet}{\raisebox{0.05cm}{$\scriptscriptstyle \bullet$}}
\newcommand{\verk}{\mbox{\scriptsize $\,\circ\,$}}
\newcommand{\halb}{\frac{1}{2}}

\newtheorem{theorem}{Theorem}[section]
\newtheorem{lemma}[theorem]{Lemma}
\newtheorem{sublemma}[theorem]{Sublemma}
\newtheorem{prop}[theorem]{Proposition}
\newtheorem{defn}[theorem]{Definition}
\newtheorem{defthm}[theorem]{Definition--Theorem}
\newtheorem{cor}[theorem]{Corollary}
\newtheorem{example}[theorem]{Example}

\newenvironment{remarknn}{\noindent {\bf Remark}}{}

\newtheorem{punkt}[theorem]{$\!\!$}
\newenvironment{proofof}{\noindent {\bf Proof of}}{\mbox{}\hfill$\Box$}
\newenvironment{proof}{\noindent {\bf Proof}}{\mbox{}\hfill$\Box$}
\begin{document}
\begin{center}
{\Large\bf A dynamical Lefschetz trace formula for algebraic Anosov diffeomorphisms}\\[0.2cm]
by Anton Deitmar and Christopher Deninger
\end{center}

\section{Introduction}
Dynamical Lefschetz trace formulas for flows on compact manifolds were first conjectured by Guillemin \cite{G} and later but independently by Patterson \cite{P}. These formulas express the alternating sum of the traces of the induced flow on cohomology by a sum over contributions from the compact orbits.

Here the cohomology theory depends on the choice of a foliation which is respected by the flow.

This note makes a contribution to Patterson's setup. He looks at Anosov flows and considers the unstable foliation. An interesting example of this situation is provided by the geodesic flow on the sphere bundle of a cocompact quotient of $\PSL_2 (\R)$ by a lattice. Using representation theory, Patterson was able to define a trace on the infinite dimensional foliation cohomologies. The dynamical Lefschetz trace formula was then shown to be a consequence of the Selberg trace formula. Incidentally this example also fits into Guillemin's framework who treated it in \cite{G}.

Much more work in the context of more general symmetric spaces was done by Bunke, Deitmar, Juhl, Olbrich and Schubert. We refer the reader to the book of Juhl \cite{J} for a comprehensive overview.

Among all Anosov flows there are the ones with an integrable complement or equivalently, those where the sum of the stable and the unstable bundle is integrable.
One gets examples by suspending Anosov diffeomorphisms and conjecturally all Anosov flows with an integrable complement should arise in this way up to a time change by a constant factor, \cite{Pl} \S\,3. 

Granting this conjecture the dynamical Lefschetz trace formula for such Anosov flows can be reduced to the following problem:

Let $f$ be an Anosov diffeomorphism of a compact manifold $X$ whose unstable foliation $\Fh^u$ is smooth. Establish a dynamical Lefschetz trace formula of the following form:
\begin{equation}
  \label{eq:1}
  \sum_p (-1)^p \Tr (f^* \tei H^p_{\Fh^u} (X)) = \sum_{x \atop f(x) = x} \varepsilon_x \, |\ddet (1- T_x f \tei T^s_x)|^{-1} \; ,
\end{equation}
where
\[
\varepsilon_x = \sgn \det (1 - T_x f \tei T_x^u) \; .
\]
Here for a fixed point $x , T^s_x$ resp. $T^u_x$ is the subspace of $T_x X$ corresponding to the eigenvalues of $T_x f$ of absolute value smaller resp. greater than one.

Note that since $f$ is supposed to be Anosov, in particular the eigenvalue $\lambda = 1$ does not occur on $T_x X$. Hence all determinants are non-zero. The definition of $H^p_{\Fh^u} (X)$ can be found in \S\,2 below.

Formula (\ref{eq:1}) is a dynamical Lefschetz trace formula for the discrete dynamical system determined by $f$.

The analytic difficulties regarding dynamical Lefschetz trace formulas for flows persist in formula (\ref{eq:1}) for diffeomorphisms. The main problem being that the complex defining $H^{\hullet}_{\Fh^u} (X)$ is not elliptic but only elliptic along the unstable leaves. In particular we cannot appeal to the work of Atiyah and Bott on Lefschetz trace formulas for endomorphisms of elliptic complexes. 

According to \cite{M} \S\,4 the following classification theorem for Anosov diffeomorphisms is expected. After lifting to a finite covering any Anosov diffeomorphism should become topologically conjugate to an algebraic Anosov diffeomorphism $f$. The latter come about as follows. Let $G$ be a simply connected real nilpotent Lie group with a cocompact lattice $\Gamma$ and consider an automorphism $f$ of $G$ with $f (\Gamma) = \Gamma$. Assume that the induced Lie algebra automorphism $f_*$ of $\eg$ does not have any eigenvalue of absolute value one. Then the induced diffeomorphism $f$ of the compact nilmanifold $X = \Gamma \ohne G$ is Anosov, and the Anosov diffeomorphisms obtained in this way are called ``algebraic''. 

If we only demand that $f$ is an endomorphism of $G$ with $f (\Gamma) \subset \Gamma$ and such that $f_*$ does not have $1$ as an eigenvalue then the induced map $f$ of $X$ will be called a generalized algebraic Anosov map.

In the present note we first prove formula (\ref{eq:1}) for algebraic Anosov diffeomorphisms. In this case by a subtle argument we can calculate the reduced leafwise cohomology i.e. the maximal Hausdorff quotient of $H^{\hullet}_{\Fh^u} (X)$. It turns out to be finite dimensional so that the definition of the trace of $f^*$ does not present a problem. 

Next, using basic aspects of the representation theory of $\eg$, we prove formula (\ref{eq:1}) also for generalized algebraic Anosov maps. In this case, even the reduced leafwise cohomology is infinite dimensional in general. Our definition of the trace of $f^*$ then depends in an essential way on representation theory. This seems to be the {\it only} method to define these traces. Namely the simplest non-trivial case of a generalized Anosov map $f$ which is not Anosov is multiplication by $-1$ on $\R / \Z$. We are not aware of any purely functional-analytic definition of a trace which applies to the map $f^*$ on $C^{\infty} (\R / \Z)$.

We would like to thank the referees whose comments have led to some improvements of the paper.

\section{Background on leafwise cohomology}
For a smooth foliation $\Fh$ on a closed smooth manifold $X$ consider the de~Rham complex of smooth real valued forms along the leaves:
\begin{equation}
  \label{eq:2}
  (\Gamma (X , \Lambda^{\hullet} T^* \Fh) , d_{\Fh}) \; .
\end{equation}
It is a complex of nuclear Fr\'echet spaces.

Here the exterior differential $d_{\Fh}$ is defined by the same formula as the ordinary differential $d$. Integrability of $T \Fh$ translates into the condition that $d^2_{\Fh} = 0$.

The cohomology groups
\begin{equation}
  \label{eq:3}
  H^{\hullet}_{\Fh} (X) := H^{\hullet} (\Gamma (X, \Lambda^{\hullet} T^* \Fh) , d_{\Fh})
\end{equation}
are called the leafwise cohomology groups of the foliated manifold $X$. If $\Rh$ is the sheaf of smooth real valued functions on $X$ which are locally constant on the leaves, we have
\begin{equation}
  \label{eq:4}
  H^{\hullet}_{\Fh} (X) = H^{\hullet} (X ,\Rh) \; .
\end{equation}
The reduced leafwise cohomology is the maximal Hausdorff-quotient $\oH^{\hullet}_{\Fh} (X)$ of $H^{\hullet}_{\Fh} (X)$. It is the nuclear Fr\'echet space obtained by dividing $\Ker d_{\Fh}$ by the closure of $\Imm d_{\Fh}$. In general even the reduced cohomology is infinite dimensional, compare example \ref{t36} for example.

Given another foliated manifold $X' , \Fh'$, any map $f : X \to X'$ which maps leaves of $\Fh$ into leaves of $\Fh'$ induces a continuous pullback map on cohomology
\[
f^* : H^{\hullet}_{\Fh'} (X') \longrightarrow H^{\hullet}_{\Fh} (X) \; .
\]
These maps pass to the reduced versions of leafwise cohomology.

Consider now a smooth compact manifold $X = \Gamma \ohne G$ obtained as the quotient of a connected real Lie group $G$ by a cocompact lattice $\Gamma \subset G$. For an immersed subgroup $P \subset G$ with Lie algebra $\ep \subset \eg$ consider the foliation $\Fh$ of $X$ with leaves $\Gamma g P$ for $g \in G$. Then we have a topological isomorphism with Lie algebra cohomology
\begin{equation}
  \label{eq:5}
  H^{\hullet}_{\Fh} (X) \silo H^{\hullet} (\ep , C^{\infty} (X)) \; .
\end{equation}
Indeed let $L = \ep \otimes C^{\infty} (X)$ be the Lie algebra of vector fields tangent to the leaves. Then we have:
\[
\Gamma (X , \Lambda^{\hullet} T^* \Fh) = \Hom_{C^{\infty} (X)} (\Lambda^{\hullet}_{C^{\infty}(X)} (L) , C^{\infty} (X)) \; .
\]
This is an isomorphisms of complexes if the right hand side is equipped with the differential (3A.5) in \cite{MS} Ch. III, Appendix. Formula (\ref{eq:5}) follows.

We now decompose the reduced leafwise cohomology using representation theory. The existence of the cocompact lattice $\Gamma$ implies that Haar measure on $G$ is unimodular. Let $\mu$ denote the unique $G$-invariant measure on $X$ of volume one. Let $R$ be the ``regular'' unitary representation of $G$ on the complex Hilbert space $L^2 (X) = L^2 (X , \mu)$ by right translation:
\[
(R_g f) (\Gamma g') = f (\Gamma g'g) \; .
\]
Since $X$ is compact, the representation $R$ decomposes into isotypic components $H (\pi)$ of finite multiplicity:
\begin{equation}
  \label{eq:6}
  L^2 (X) = \widehat{\bigoplus_{\pi\in \hat{G}}} H (\pi) \; .
\end{equation}
Here $\hat{G}$ is the unitary dual of $G$, i.e. the set of equivalence classes of irreducible unitary representations of $G$. Each isotypic component is isomorphic to $M_{\pi} \otimes V_{\pi}$ where $(\rho_{\pi} , V_{\pi})$ is a fixed representation in the class $\pi$ and where $M_{\pi}$ is the finite dimensional multiplicity space, i.e.
\[
M_{\pi} = \Hom_G (V_{\pi} , L^2 (X)) \; .
\]
We can now analyse the reduced leafwise cohomology using representation theory. Denote by $\oH^p (\ep , H (\pi)^{\infty})$ the maximal Hausdorff quotient of \\$H^p (\ep , H (\pi)^{\infty})$, where $H (\pi)^{\infty}$ denotes the space of smooth vectors in $H (\pi)$. 

\begin{lemma}
  \label{t21}
The natural maps 
\[
\oH^p (\ep , H (\pi)^{\infty}) \longrightarrow \oH^p_{\Fh} (X, \C) = \oH^p_{\Fh} (X) \otimes \C
\]
are split injections and hence have closed image. In this way the algebraic direct sum
\[
\bigoplus_{\pi} \oH^p (\ep , H (\pi)^{\infty})
\]
becomes a dense subspace of
\[
\oH^p_{\Fh} (X , \C) = \oH^p_{\Fh} (X) \otimes \C \; .
\]
Every element of $\oH^p_{\Fh} (X , \C)$ can be written as an unconditionally convergent series $\sum_{\pi} h_{\pi}$ with $h_{\pi} \in \oH^p (\ep , H (\pi)^{\infty})$.
\end{lemma}

The proof is given below. It is based on an auxiliary result which requires the following notions:

\begin{defn}\label{t22_neu}
  A sequence $(E_i)_{i \ge 1}$ of non-trivial closed linear subspaces of a Fr\'echet space $E$ is called a Schauder basis of subspaces of $E$ if the following conditions hold:\\
a) Every element $v \in E$ can be written uniquely as a convergent series $v = \sum^{\infty}_{i=1} v_i$ with $v_i \in E_i$.\\
b) The projections $E \to E_i , v \mapsto v_i$ resulting from a) are continuous.\\
If the series in a) converge unconditionally the basis $(E_i)$ is called unconditional.
\end{defn}
\begin{sublemma}
  \label{t22}
Let $C^{\hullet}$ be a complex of Fr\'echet spaces with continuous differentials. Consider closed subcomplexes $C^{\hullet}_i \subset C^{\hullet}$ for $i \ge 1$ such that for every $p$ the sequence $(C^p_i)_{i \ge 1}$ is an unconditional Schauder basis of subspaces of $C^p$.
Then the natural maps 
\[
\oH^p (C^{\hullet}_i) \longrightarrow \oH^p (C^{\hullet})
\]
are split injective and hence have closed image. They induce an inclusion with dense image
\[
\bigoplus_i \oH^p (C^{\hullet}_i) \subset \oH^p (C^{\hullet}) \; .
\]
Every element of $\oH^p (C^{\hullet})$ can be written as an unconditionally convergent series $\sum_i h_i$ with $h_i \in \oH^p (C^{\hullet}_i)$. 
\end{sublemma}

\begin{proof}
  For any $p$ the assumptions give continuous projections
\[
\beta_i : C^{p-1} \longrightarrow C^{p-1}_i
\]
such that for all $c$ in $C^{p-1}$ we have
\[
c = \sum_i \beta_i (c) \; .
\]
Since the differential is continuous we get
\[
dc = \sum_i d\beta_i (c) \; .
\]
On the other hand we have
\[
dc = \sum_i \beta_i (dc) \; .
\]
The representation of an element in $C^p$ as a convergent series of elements in the $C^p_i$ is unique. Hence we find
\[
d\beta_i (c) = \beta_i (dc)
\]
for all $c$ i.e. $d \verk \beta_i = \beta_i \verk d$. From this and the continuity of $\beta_i$ it follows that
\begin{equation} \label{eq:7_neu}
\beta_i (\overline{d C^{p-1}}) \subset \overline{\beta_i (d C^{p-1})} = \overline{d \beta_i (C^{p-1})} = \overline{d C^{p-1}_i} \; .
\end{equation}
The inclusion $\kappa_i : C^{\hullet}_i \hookrightarrow C^{\hullet}$ of complexes induces a continuous map on the $p$-th cohomology and hence a continuous map
\[
\okappa_i : \oH^p (C^{\hullet}_i) \longrightarrow \oH^p (C^{\hullet}) \; .
\]
Similarly the continuous surjective map of complexes 
\[
\beta_i : C^{\hullet} \longrightarrow C^{\hullet}_i
\]
leads to a continuous map
\[
\obeta_i : \oH^p (C^{\hullet}) \longrightarrow \oH^p (C^{\hullet}_i) \; .
\]
We have $\beta_i \verk \kappa_i = \id$ and hence $\obeta_i \verk \okappa_i = \id$. Therefore $\okappa_i$ is a split injection. Its image is closed since it is also the kernel of the continuous projector $\pi_i = \id - \okappa_i \verk \obeta_i$ on $\oH^p (C^{\hullet})$. 

We now show that the induced map
\[
\bigoplus_i \oH^p (C^{\hullet}_i) \longrightarrow \oH^p (C^{\hullet})
\]
is injective. So assume that the element $\sum_i [c_i]$ goes to zero in $\oH^p (C^{\hullet})$. This means that $\sum_i c_i$ lies in $\overline{dC^{p-1}}$. Here $c_i \in C^p_i , dc_i = 0$ and almost all $c_i$ are zero.  Applying $\beta_i$ we see that $c_i$ lies in $\beta_i (\overline{dC^{p-1}})$ and by (\ref{eq:7_neu}) therefore in $\overline{dC^{p-1}_i}$. Thus $[c_i] = 0$ for all $i$.

Next consider any element $[c]$ in $\oH^p (C^{\hullet})$. We may write $c$ as an unconditionally convergent series
\[
c = \sum_i c_i \quad \mbox{with} \; c_i = \beta_i (c) \in C^p_i \; .
\]
We have $0 = \beta_i (dc) = d \beta_i (c) = dc_i$. Hence $[c_i] \in \oH^p (C^{\hullet})$ lies in the image of $\oH^p (C^{\hullet}_i)$ and we have
\[
[c] = \sum_i [c_i]
\]
in $\oH^p (C^{\hullet})$ the series being unconditionally convergent. In particular the algebraic direct sum of the spaces $\oH^p (C^{\hullet}_i)$ is dense in $\oH^p (C^{\hullet})$.
\end{proof}

\begin{proofof} {\bf \ref{t21}} 
Using the Sobolev embedding theorem one sees that the non-zero subspaces $H (\pi)^{\infty} \subset C^{\infty} (X) = L^2 (X)^{\infty}$ form an unconditional Schauder basis of subspaces of $C^{\infty} (X)$. It follows that the subcomplexes $C^{\hullet} (\eg , H (\pi)^{\infty})$ of $C^{\hullet} (\eg , C^{\infty} (X))$ satisfy the assumptions of the sublemma. This implies lemma \ref{t21}. 
\end{proofof}


\section{$\Gamma$-acceptable subgroups and leafwise cohomology}

In this section we first prove that the reduced leafwise cohomology of certain homogenous foliations on nilmanifolds is isomorphic to finite dimensional Lie algebra cohomology. This result applies to the stable and unstable foliations of an algebraic Anosov diffeomorphism. It also applies to the foliation defined by a normal subgroup where it gives a result of \'Alvarez L\'opez and Hector.

Let $G$ be a simply connected nilpotent Lie group with a cocompact lattice $\Gamma \subset G$. Let
\begin{equation}
  \label{eq:8_neu}
  G = C_0 \supsetneqq C_1 \supsetneqq \ldots \supsetneqq C_k \supsetneqq C_{k+1} = 1
\end{equation}
be the descending central series of $G$ defined by $C_{j+1} = [G , C_j]$ and $C_0 = G$. From the definition we see:

\begin{punkt}\label{t31_neu}
  \rm The group $C_j / C_{j+1}$ is central in $G / C_{j+1}$ for every $0 \le j \le k$. In particular $C_k$ is central in $G$.
\end{punkt}
We will need the following basic results of Malcev \cite{Ma}:
\begin{punkt}\label{t32_neu}
  \rm The group $\Gamma_j = \Gamma \cap C_j$ is a cocompact lattice in $C_j$ for every $j$.
\end{punkt}
\begin{punkt}\label{t33_neu}
  \rm The image $\Gamma (j)$ of $\Gamma$ in $G / C_j$ is a cocompact lattice for every $j$.
\end{punkt}
It follows from \ref{t32_neu} and \ref{t33_neu} that the image of $\Gamma_j$ in the abelian Lie algebra $C_j / C_{j+1}$ is a cocompact lattice. The quotient is a real torus
\[
\Delta_j = C_j / C_{j+1} \Gamma_j = \Gamma_j \ohne C_j / C_{j+1} \; .
\]
\begin{defn}
  \label{t34_neu}
A closed connected subgroup $P$ of $G$ is called $\Gamma$-acceptable if the image of $P_j = P \cap C_j$ in $\Delta_j$ is dense for every $0 \le j \le k$.
\end{defn}
In \ref{t37_neu} below we will see that the stable and unstable groups of an algebraic Anosov diffeomorphism of $X = \Gamma \ohne G$ are $\Gamma$-acceptable. Moreover, normal subgroups of $G$ whose orbits in $X$ are dense are $\Gamma$-acceptable by proposition \ref{t311_neu}.

Let $\Fh$ be the foliation of $X$ by the orbits of the right $P$-operation on $X$. The embedding of $\R$ into $C^{\infty} (X)$ induces a map
\begin{equation}
  \label{eq:9_neu}
  H^{\hullet} (\ep , \R) \longrightarrow \oH^{\hullet} (\ep , C^{\infty} (X)) \overset{(\ref{eq:5})}{\cong} \oH^{\hullet}_{\Fh} (X) \; .
\end{equation}
\begin{theorem}
  \label{t35_neu}
If $P$ is $\Gamma$-acceptable this is an isomorphism:
\[
H^{\hullet} (\ep , \R) \cong \oH^{\hullet}_{\Fh} (X) \; .
\]
\end{theorem}
For the proof we require the following lemma.

\begin{lemma}
  \label{t36_neu}
Let $\ep$ be a finite dimensional real Lie algebra and let
\[
0 = \ea_0 \subset \ldots \subset \ea_n \subset \ep
\]
be a chain of normal subalgebras of $\ep$ such that $\ea_j / \ea_{j-1}$ is central in $\ep / \ea_{j-1}$ for every $1 \le j \le n$. Let $M$ be a complex vector space on which $\ep$ acts $\C$-linearly. We assume that $\ea_{n-1}$ acts trivially and that $\ea_n$ acts by a nontrivial character $\alpha$. Then $H^{\hullet} (\ep , M) = 0$.
\end{lemma}

\begin{proof}
  The Hochschild--Serre spectral sequence
\[
H^p (\ep / \ea_n , H^q (\ea_n , M)) \Longrightarrow H^{p+q} (\ep , M)
\]
shows that we may assume that $\ep = \ea_n$. Clearly we may also assume that the inclusions in the chain
\[
0 = \ea_0 \subset \ldots \subset \ea_n = \ep
\]
are strict. Since $\ea_n = \ep$ acts non-trivially we have $\ep \neq 0$ and hence $n \ge 1$. The assertion is proved by induction on $n$. For $n = 1$ the Lie algebra $\ep$ is abelian. Hence it operates on $M$-valued cocycles and therefore on $H^{\hullet} (\ep , M)$ by multiplication with $\alpha$ c.f. \cite{HS} p. 591 below. On the other hand the operation of $\eg$ on $H^{\hullet} (\eg , N)$ is trivial for any Lie algebra $\eg$ and any $\eg$-module $N$. Since $\alpha$ was non-trivial, it follows that $H^{\hullet} (\ep , M) = 0$. 

Now fix some $n \ge 2$ and assume the assertion proven for $n - 1$. Consider the spectral sequence 
\begin{equation}
  \label{eq:10_neu}
  H^p (\ep / \ea_1 , H^q (\ea_1 , M)) \Longrightarrow H^{p+q} (\ep , M) \; .
\end{equation}
Since $\ea_1$ is central in $\ep = \ea_n$ it follows that $\ep$ acts by multiplication with $\alpha$ on $M$-valued cocycles of $\ea_1$ and hence on $H^q (\ea_1 , M)$ c.f. \cite{HS} p. 591. Setting $M' = H^q (\ea_1 , M)$ and $\ea'_i = \ea_{i+1} / \ea_0 , \ep' = \ep / \ea_1$ we get a chain
\[
0 = \ea'_0 \subset \ldots \subset \ea'_{n-1} = \ep' \; .
\]
Applying the induction hypotheses to the $\ep'$-action on $M'$ gives
\[
H^{\hullet} (\ep / \ea_1 , H^q (\ea_1 , M)) = H^{\hullet} (\ep' , M') = 0 \; .
\]
Using the spectral sequence (\ref{eq:10_neu}) we conclude that $H^{\hullet} (\ep , M) = 0$.
\end{proof}

\begin{proofof}
  {\bf \ref{t35_neu}} It is enough to show that the natural inclusion $\C \hookrightarrow \C^{\infty} (X) = C^{\infty} (X , \C)$ induces an isomorphism
  \begin{equation}
    \label{eq:11_neu}
    H^{\hullet} (\ep , \C) \silo \oH^{\hullet} (\ep , \C^{\infty} (X)) \; .
  \end{equation}
According to \ref{t33_neu} the image $\Gamma (j)$ of $\Gamma$ in $G (j) = G / C_j$ is a lattice for every $0 \le j \le k+1$. The natural projection of right $P$-spaces
\[
X \longrightarrow X_j = \Gamma (j) \ohne G (j) = \Gamma \cdot C_j \ohne G
\]
induces a $\ep$-equivariant injection
\begin{equation}
  \label{eq:12_neu}
  \C^{\infty} (X_j) \hookrightarrow \C^{\infty} (X) \; .
\end{equation}
We will show by descending induction on $j$ that (\ref{eq:12_neu}) induces isomorphisms of Fr\'echet spaces for all $0 \le j \le k +1$
\begin{equation}
  \label{eq:13_neu}
  \oH^{\hullet} (\ep , \C^{\infty} (X_j)) \silo \oH^{\hullet} (\ep , \C^{\infty} (X)) \; .
\end{equation}
For $j = k+1$ this is clear. For the induction step fix $0 \le m \le k$ and assume that (\ref{eq:13_neu}) is an isomorphism for $j = m+1$. Now consider the torus
\[
\Delta = \Delta_m = C_m / C_{m+1} \Gamma_m = \Gamma_m \ohne C_m / C_{m+1} \; .
\]
The right $G$-action on $X_{m+1}$ induces a $\Delta$-action on $X_{m+1}$ since the elements of $C_m$ are central $\mod C_{m+1}$. Fourier theory implies that the $\chi$-isotypical components
\[
\C^{\infty} (X_{m+1}) (\chi)
\]
of the induced $\Delta$-action on $\C^{\infty} (X_{m+1})$ form an unconditional Schauder basis of subspaces in the sense of definition \ref{t22_neu}. Here $\chi$ runs over the characters of $\Delta$. It follows that the subcomplexes
\[
C^{\hullet} (\ep , \C^{\infty} (X_{m+1}) (\chi))
\]
of $C^{\infty} (\ep , \C^{\infty} (X_{m+1}))$ satisfy the assumptions of sublemma \ref{t22}. Hence there is a natural injection
\begin{equation}
  \label{eq:14_neu}
  \bigoplus_{\chi} \oH^{\hullet} (\ep , \C^{\infty} (X_{m+1}) (\chi)) \subset \oH^{\hullet} (\ep , \C^{\infty} (X_{m+1}))
\end{equation}
with dense image. Moreover the spaces $\oH^{\hullet} (\ep , \C^{\infty} (X_{m+1}) (\chi))$ are closed in $\oH^{\hullet} (\ep , \C^{\infty} (X_{m+1}))$. We now apply lemma \ref{t36_neu} to $\ep = \Lie P$ and the chain
\[
0 = \ea_0 \subset \ldots \subset \ea_n \subset \ep
\]
defined by $\ea_j = \Lie (P \cap C_{k+1-j})$. Here $n = k+1 -m$. For $M$ we take $M = \C^{\infty} (X_{m+1}) (\chi)$ where $\chi$ is any character of $\Delta$ which is non-trivial on $P_m = P \cap C_m$. Then $\ea_n = \Lie P_m$ acts by the non-trivial character $\alpha = \chi_*$ on $M$ and $\ea_{n-1} = \Lie P_{m+1}$ acts trivially. It follows from the lemma that only $\chi$'s with $\chi (P_m) = 1$ can contribute to the direct sum in (\ref{eq:14_neu}). 

Now, by assumption $P$ is $\Gamma$-acceptable. Hence the image of $P_m$ in $\Delta$ is dense. The condition $\chi (P_m) = 1$ therefore implies that $\chi$ is the trivial character $\chi = \1$. Thus $\oH^{\hullet} (\ep , \C^{\infty} (X_{m+1}) (\1))$ is a dense subspace of $\oH^{\hullet} (\ep , \C^{\infty} (X_{m+1}))$. Since it is also closed we get:
\begin{equation}
  \label{eq:15_neu}
  \oH^{\hullet} (\ep , \C^{\infty} (X_m)) = \oH^{\hullet} (\ep , \C^{\infty} (X_{m+1})) \; .
\end{equation}
Here we took into account that
\[
\C^{\infty} (X_{m+1}) (\1) = \C^{\infty} (X_{m+1} / \Delta) = \C^{\infty} (X_m) \; .
\]
By the induction hypotheses we have an isomorphism
\[
\oH^{\hullet} (\ep , \C^{\infty} (X_{m+1})) \silo \oH^{\hullet} (\ep , \C^{\infty} (X)) \; .
\]
Combining it with (\ref{eq:15_neu}) we get the isomorphism (\ref{eq:13_neu}) for $j = m$. Note that a continuous algebraic isomorphism of Fr\'echet spaces is necessarily a topological isomorphism by the open mapping theorem. Thus we have completed the induction step. The case $j = 0$ of (\ref{eq:13_neu}) now give the assertion of the theorem.
\end{proofof}

\begin{punkt}
  \label{t37_neu} \rm
Let $f$ be an algebraic Anosov diffeomorphism of a nilmanifold $X = \Gamma \ohne G$ as defined in the introduction. It is induced by an automorphism of $G$ which we also denote by $f$. Let
$\eg^s \subset\eg$ be the subspace of all $v\in\eg$ for which
$f_*^n(v)$ tends to zero as $n\rightarrow\infty$. This is the {\it
stable} subspace of $\eg$. The {\it unstable} subspace $\eg^u$ of $\eg$ is the space of all $v$ of $\eg$ such that $f^n_* (v) \to 0$ for $n \to - \infty$. Since $f$ is a group
automorphism it follows that $f_*$ is a Lie algebra
automorphism and so $\eg^s$ and $\eg^u$ are Lie subalgebras of $\eg$. Let
$G^s$ and $G^u$ denote the corresponding connected subgroups of $G$.
\end{punkt}

\begin{prop}
  \label{t38_neu}
In the situation of \ref{t37_neu} the groups $P = G^s$ and $P = G^u$ are $\Gamma$-acceptable.
\end{prop}

\begin{proof}
The stable subgroup of $f^{-1}$ is the unstable subgroup of $f$. Hence it is enough to treat the case $P = G^s$.  We have to show that for $0 \le j \le k$ the image of $P_j$ is dense in $\Delta_j = \Gamma_j \ohne G_j$ where $G_j = C_j / C_{j+1}$ is abelian. We have $P_j = G^s_j$ with respect to the induced automorphism $f_j$ of $G_j$. Hence it suffices to show that for abelian $G$ the image of $P = G^s$ is dense in $\Gamma \ohne G$ or equivalently that the image of $\Gamma$ is dense in $G / G^s$.

So suppose $G$ is abelian. Then the exponential map is an
isomorphism of the Lie algebra of $G$ with $G$, so $G$ is
a finite dimensional $\R$-vector space. We therefore
rephrase the claim as follows. Let $V$ be a non-zero finite
dimensional real vector space and let $f$ be an
automorphism of $V$ with no complex eigenvalue of absolute
value $1$. Let $V^s$ be the stable subspace, i.e.,
$$
V^s\=\{ v\in V\mid f^n(v)\rightarrow 0\ {\rm as}\ n\rightarrow\infty\}.
$$
Let $\Gamma \subset V$ be a lattice
with $f(\Gamma)= \Gamma$. Then we have to show that the image $\Gamma_W$ of
$\Gamma$ is dense in $W = V/V^s$. This will follow if for every $\varepsilon$-neighborhood $U$ of zero in $W$ the set $U \cap \Gamma_W$ contains a basis of $W$. Namely the lattice generated by such a basis is contained in $\Gamma_W$ and has a fundamental domain in $U$. Thus every given point in $W$ is $\varepsilon$-close to some point of this lattice and hence of $\Gamma_W$. \\
So, let $U$ be a neighborhood of zero in $W$. On $W$ all eigenvalues of $f$
are of absolute value $>1$, therefore $f^{-1}$ is
contracting there. Let $v_1,\dots , v_r$ be a basis of $W$
contained in $\Gamma_W$. Then the vectors $f^{-n}(v_1),\dots,f^{-n}(v_r)$
form a basis for every $n\in\N$ and this basis shrinks
to zero as $n\rightarrow\infty$, so the claim follows.
\end{proof}

\begin{remarknn}
  Indirectly the proof shows that $V^s \neq 0$. This is clear a priori since on the one hand $| \det f | = 1$ because $f$ respects a lattice and on the other hand $f$ does not have eigenvalues of absolute value $1$. 
\end{remarknn}
\begin{cor}
  \label{t39_neu}
In the situation of \ref{t37_neu} we have natural isomorphisms 
\[
\oH^{\hullet}_{\Fh^s} (X) \cong H^{\hullet} (\eg^s , \R) \quad \mbox{and} \quad \oH^{\hullet}_{\Fh^u} (X) \cong H^{\hullet} (\eg^u , \R) \; .
\]
In particular these reduced leafwise cohomologies are finite dimensional.
\end{cor}

\begin{proof}
  This follows from proposition \ref{t38_neu} and theorem \ref{t35_neu} using that $\Fh^s$ and $\Fh^u$ are the foliations by the orbits of the $G^s$-resp. $G^u$-action on $X$.
\end{proof}

\begin{remarknn}
  We conjecture that for {\it any} Anosov diffeomorphism of a compact manifold with smooth stable and unstable foliations the reduced leafwise cohomologies with respect to these foliations are finite dimensional.
\end{remarknn}

\begin{punkt} \rm
  \label{t310_neu}
In this final subsection we consider normal subgroups in simply connected real Lie groups. The following result holds:
\end{punkt}

\begin{prop}
  \label{t311_neu}
A normal subgroup $P$ with $\Gamma P$ dense in $G$ is
$\Gamma$-acceptable.
\end{prop}

\begin{remarknn}
  The proposition and theorem \ref{t35_neu} imply the isomorphism
\[
\oH^{\hullet}_{\Fh} (X) \cong H^{\hullet} (\ep , \R) \; .
\]
This result is due originally to \'Alvarez L\'opez and Hector \cite{AH} Theorem 2.12 with a different proof.
\end{remarknn}

\begin{proof} 
It suffices to show that $\Gamma_j P_j$ is dense
in $C_j$. We use induction on $j$. For $j=0$ this is the assumption. Assume the claim proven for $j$. Since
$C_{j+1}=[G,C_j]$ we pick $x\in G$ and $y\in C_j$. Let
$\gamma_n p_n$ be a sequence converging to $x$, where
$\gamma_n\in \Gamma$ and $p_n\in P$. Likewise let $\beta_n l_n$ be
a sequence with $\beta_n \in \Gamma_j$ and $l_n \in P_j$ converging to $y$. Then the
commutator $[\gamma_n p_n,\beta_n l_n]$ converges to $[x,y]$. If
we can show that the former commutator lies in
$\Gamma_{j+1} P_{j+1}$ then we are done since $P_{j+1}$ is
normal in $G$ so that $\Gamma_{j+1} P_{j+1}$ is a subgroup. Now the claim will follow from the
following group theoretical lemma.
\end{proof}

\begin{lemma}
Let $G$ be a group and let $A,B$ be subgroups. Suppose
that $B$ is normal. Let $C_j$ for $j=0,1,2,\dots$ denote
the lower central series. Let $A_j=A\cap C_j$ and
$B_j=B\cap C_j$. Then
$$
[AB,A_jB_j]\ \subset\ A_{j+1}B_{j+1}.
$$
\end{lemma}

\begin{proof} Let $a\in A, b\in B, \alpha\in A_j, \beta\in B_j$. For two
group elements $x,y$ we write $x^y$ for $y^{-1} x y$. Then
\begin{eqnarray*}
[a b,\alpha\beta] &=& ab\alpha\beta b^{-1} a^{-1} \beta^{-1}\alpha^{-1}\\
   &=& a\alpha b^\alpha \beta b^{-1} a^{-1} \beta^{-1} \alpha^{-1}\\
   &=& [a,\alpha] b^{\alpha a^{-1} \alpha^{-1}} \beta^{a^{-1}\alpha^{-1}} (b^{a^{-1} \alpha^{-1}})^{-1} (\beta^{\alpha^{-1}})^{-1} \\
   &=& [a,\alpha]\left[b^{\alpha a^{-1} \alpha^{-1}},\beta^{a^{-1}\alpha^{-1}}\right] \beta^{a^{-1}\alpha^{-1}} b^{\alpha a^{-1} \alpha^{-1}} (b^{a^{-1} \alpha^{-1}})^{-1} (\beta^{\alpha^{-1}})^{-1} \; .
\end{eqnarray*}
The two commutators are in $A_{j+1}$ resp. $B_{j+1}$ alright
so we only have to care about the rest,
$$
\beta^{a^{-1}\alpha^{-1}} b^{\alpha a^{-1}\alpha^{-1}}
(b^{-1})^{a^{-1}\alpha^{-1}} (\beta^{-1})^{\alpha^{-1}} \; .
$$
In the same fashion as before we can move the first
$\beta$-term past the two $b$-terms producing  commutators in
$B_{j+1}$. It remains to consider
$$
b^{\alpha a^{-1}\alpha^{-1}} (b^{-1})^{a^{-1} \alpha^{-1}}
 \beta^{a^{-1}\alpha^{-1}}  (\beta^{-1})^{\alpha^{-1}} \= \left(\alpha^{-1}
 b\alpha b^{-1}\right)^{a^{-1} \alpha^{-1}}\left(a\beta
 a^{-1}\beta^{-1}\right)^{\alpha^{-1}}.
$$
The commutator $\alpha^{-1} b\alpha b^{-1}=[\alpha^{-1},b]$ lies in
$[A_j,B]\subset C_{j+1}$. But on the other hand, since $B$
is normal, this commutator lies in $B$, therefore in
$B_{j+1}$. The latter also is normal, hence the first
factor lies in $B_{j+1}$. The second commutator $a\beta
a^{-1}\beta^{-1}=[a,\beta]$ lies in $[A,B_j]$ and likewise
turns out to lie in $B_{j+1}$. 
\end{proof}


\section{The dynamical Lefschetz trace formula for algebraic Anosov diffeomorphisms}

Consider an algebraic Anosov diffeomorphism of a nilmanifold $X = \Gamma \ohne G$ induced by an automorphism $f$ of $G$. According to corollary \ref{t39_neu} the reduced leafwise cohomology $\oH^p_{\Fh^u} (X)$ is finite dimensional. Thus the trace of the induced automorphism $f^*$ on this cohomology is well defined.

\begin{theorem}
  \label{t41_neu}
The following formula holds:
\[
\sum_p (-1)^p \Tr (f^* \tei \oH^p_{\Fh^u} (X)) = \sum_{x \atop f (x) = x} \varepsilon_x \, |\det (1- T_x f \tei T^s_x)|^{-1} \; .
\]
Here
\[
\varepsilon_x = \sgn \det (1 - T_x f \tei T^u_x) \; .
\]
\end{theorem}

\begin{proof}
Set $\ep = \eg^u$ and $\Fh = \Fh^u$. Then we have:
  \begin{eqnarray}
    \label{eq:9}
    \sum_p (-1)^p \Tr (f^* \tei \oH^p_{\Fh} (X)) & = & \sum_p (-1)^p \Tr (f^* \tei H^p (\ep , \R)) \nonumber \\
& = & \sum_p (-1)^p \Tr (\Lambda^p \check{f}_* \tei \Lambda^p \check{\ep}) \nonumber \\
& = & \det (1- f_* \tei \ep) \nonumber \\
& = & \det (1 - f_* \tei \eg) \det (1 - f_* \tei \eg / \ep)^{-1} \; .
  \end{eqnarray}
Now by a theorem of Nomizu \cite{N} it is known that
\[
H^p (\eg , \R) = H^p (X , \R)
\]
for all $p$. Therefore by the same calculation as before:
\begin{eqnarray*}
  \det (1 - f_* \tei \eg) & = & \sum_p (-1)^p \Tr (f^* \tei H^p (\eg , \R)) \\
& = & \sum_p (-1)^p \Tr (f^* \tei H^p (X , \R)) \; .
\end{eqnarray*}
As $f$ has only non-degenerate fixed points (see below), the ordinary Lefschetz trace formula now implies:
\begin{equation}
  \label{eq:10}
\det (1- f_* \tei \eg) = \sum_{x \atop f (x) = x} \sgn \det (1 - T_x f \tei T_x X) \; .
\end{equation}
For a fixed point $x = \Gamma g$ of $f$ we have $\Gamma = \Gamma g f(g)^{-1}$. The commutative diagram:
\[
\begin{CD}
  X @>{R_{g^{-1}}}>> X \\
@V{f}VV @VV{f}V \\
X @>>{R_{f (g)^{-1}}}> X
\end{CD}
\]
induces commutative diagrams where $\overline{e} = \Gamma \in X$:
\begin{equation}
  \label{eq:11}
  \begin{CD}
    T_x X @>{\sim}>> T_{\overline{e}} X @= \eg \\
@V{T_x f}VV @V{T_{\overline{e}} f}VV @VV{f_*}V \\
T_xX @>{\sim}>> T_{\overline{e}} X @= \eg
  \end{CD} \hspace*{1.5cm} \mbox{and} \hspace*{1.5cm}
  \begin{CD}
    T_x \Fh @>{\sim}>> \ep \\
@V{T_x f}VV @VV{f_*}V \\
T_x \Fh @>{\sim}>> \ep \; .
  \end{CD}
\end{equation}
Hence:
\[
\sum_{x \atop f (x) = x} \varepsilon_x \, |\ddet (1- T_x f \tei T_x X / T_x \Fh)|^{-1} = \frac{\sgn \det (1 - f_* \tei \ep)}{|\ddet (1 - f_* \tei \eg / \ep)|} \; \sharp \; \Fix f
\]
and (\ref{eq:10}) asserts that:
\[
\det (1 - f_* \tei \eg) = \sgn \det (1 - f_* \tei \eg) \; \sharp \; \Fix f \; .
\]
Together with equation (\ref{eq:9}) this gives the formula if we remark that $T_x X / T_x \Fh = T_x \Fh^s = T^s_x$.
\end{proof}

\section{Dynamical Lefschetz trace formulas for generalized algebraic Anosov maps}
In this section we will prove dynamical Lefschetz trace formulas for generalized algebraic Anosov maps with respect to homogenous invariant foliations. 

Thus let $f$ be a generalized algebraic Anosov map for the pair $(G , \Gamma)$ as in the introduction and let $\ep \subset \eg$ be an $f_*$-invariant sub Lie algebra of $\eg$. To $\ep$ corresponds an $f$-invariant closed subgroup $P$ of $G$. Denote by $\Fh$ the corresponding foliation on $X = \Gamma \ohne G$.

By assumption $f$ maps leaves of $\Fh$ into leaves of $\Fh$ and hence it induces an endomorphism $f^*$ of $\oH^p_{\Fh} (X)$. Under the isomorphism
\begin{equation}
  \label{eq:7}
  \oH^p_{\Fh} (X) = \oH^p (\ep , C^{\infty} (X))
\end{equation}
the operation $f^*$ corresponds to the map on Lie algebra cohomology induced by $f_*$ on $\ep$ and by pullback $f^*$ on $C^{\infty} (X)$. 

We want to define the trace of $f^*$ on the generally infinite dimensional cohomology (\ref{eq:7}). For this we use representation theory.

Let $(\rho , V)$ be any representation of $G$. The pullback representation along an endomorphism $f : G \to G$ is defined by
\[
f^* (\rho , V) = (\rho \verk f , V) \; .
\]
It is unitary if $(\rho , V)$ is. For $\pi \in \hat{G}$ we get an equivalence class $f^* \pi$ of unitary representations. If $f : G \to G$ is an isomorphism, then $f^* \pi$ is irreducible. In this case we get a bijection $f^* : \hat{G} \to \hat{G}$.

The following lemma will be crucial

\begin{lemma}
  \label{t31}
Let $f$ be an endomorphism of a simply connected nilpotent Lie group $G$ such that $f_*$ does not have $1$ as an eigenvalue. For $\pi \in \hat{G}$ assume that there is a non-zero $G$-intertwining map $T$ from $f^* (V^{\infty}_{\pi})$ to $V^{\infty}_{\pi}$. Then $\pi$ is the trivial representation.
\end{lemma}

\begin{proof}
 Let us write $(\rho , V)$ for $(\rho^{\infty}_{\pi} , V^{\infty}_{\pi})$. Let $C_j$ be the descending central series. The integer $k \ge 0$ such that $C_{k+1} = 1$ but $C_k \neq 1$ is called the length of the central series. Since $f$ is an endomorphism it preserves all the $C_j$ and in particular the subgroup $C_k$. This connected subgroup lies in the center of $G$. Since $\pi$ was supposed to be irreducible there is hence a character
\[
\chi : C_k \longrightarrow \C^*
\]
such that
\[
\rho (g) = \chi (g) \, \id_V \quad \mbox{for all} \; g \in C_k \; .
\]
Since $f$ preserves $C_k$ it follows that
\[
(\rho \verk f) (g) = (\chi \verk f) (g) \, \id_V \quad \mbox{for all} \; g \in C_k \; .
\]
The intertwining relation $\rho (g) T = T (\rho \verk f) (g)$ for all $g \in  G$ now implies that
\[
\chi (g) T = (\chi \verk f) (g) T \quad \mbox{for all} \; g \in G \; .
\]
Since $T$ was supposed to be non-trivial this gives the relation $\chi =\chi \verk f$ on $C_k$. For the derivative $\chi_* : \eZ_k = \Lie C_k \to \C$ it follows that
\[
\chi_* = \chi_* \verk f_* \quad \mbox{on} \; \eZ_k \; .
\]
By assumption $f_*$ does not have eigenvalue $1$ on $\eg$ and all the less on $\eZ_k$. Hence the transposed map $\check{f}_*$ does not have eigenvalue $1$ on $\check{\eZ}_k$. It follows that $\chi_* = 0$, so $\chi$ is trivial.

If $k = 0$ then $C_k = G$ is abelian and it follows that $\rho = \chi = 1$. Thus $\pi$ is trivial.

If $k > 0$ then $\rho$ factors over $G / C_k$ which is nilpotent of length $k -1$. The map $f$ induces an endomorphism of this smaller group whose derivative does not have the eigenvalue $1$. Also, $T$ remains an intertwining map for $G / C_k$. Thus we can iterate the argument down to the abelian case $k = 0$ to reach the same conclusion. 
\end{proof}

Recall the decomposition (\ref{eq:6}). According to lemma \ref{t21} the algebraic direct sum
\[
\Hh = \bigoplus_{\pi} \oH^p (\ep , H (\pi)^{\infty})
\]
is a dense subspace of $\oH^p_{\Fh} (X, \C)$.

\begin{defthm}
  \label{t32}
Let $f$ be a generalized algebraic Anosov map for $(G, \Gamma)$. Then we have the following assertions:\\
1) If $f : G \to G$ is surjective, hence an automorphism, then the map $f^*$ on $\oH^p_{\Fh} (X,\C)$ permutes the closed subspaces $\oH^p (\ep , H (\pi)^{\infty})$:
\[
f^* (\oH^p (\ep , H (\pi)^{\infty})) = \oH^p (\ep , H (f^* \pi)^{\infty}) \; .
\]
In particular $f^*$ leaves $\Hh$ invariant. The condition $f^* \pi = \pi$ for the subspace $\oH^p (\ep , H (\pi)^{\infty})$ to be left invariant is equivalent to $\pi$ being the trivial representation. Hence it is reasonable to define the trace of $f^*$ as follows:
\begin{eqnarray*}
  \Tr (f^* \tei \oH^p_{\Fh} (X)) & := & \Tr (f^* \tei \Hh) \\
& := & \sum_{f^* \pi = \pi} \Tr (f^* \tei \oH^p (\ep , H (\pi)^{\infty})) \\
& = & \Tr (f^* \tei H^p (\ep , \R)) \; .
\end{eqnarray*}
If $\oH^p_{\Fh} (X)$ is finite dimensional then our definition agrees with the usual trace of $f^*$ as an endomorphism of a finite dimensional space.

2) In general, if $f$ is only an endomorphism of $G$, then unless $\pi$ is the trivial representation,
\[
f^* (H (\pi)^{\infty}) \quad \mbox{is orthogonal to} \; H (\pi)^{\infty}
\]
in the Hilbert space $L^2 (X)$. Here
\[
f^* : C^{\infty} (X) \longrightarrow C^{\infty} (X)
\]
is the pullback map induced by $f : X \to X$. It follows that for non-trivial $\pi$ we have:
\[
f^* \oH^p (\ep , H (\pi)^{\infty}) \subset \overline{\bigoplus_{\pi' \neq \pi}} \oH^p (\ep , H (\pi')^{\infty}) \; ,
\]
where $\overline{\bigoplus}$ denotes the closure of the direct sum in the topology of $\oH^p_{\Fh} (X)$. Again this suggests to define the trace of $f^*$ by the formula
\[
\Tr (f^* \tei \oH^p_{\Fh} (X)) := \Tr (f^* \tei \oH^p (\ep , \R)) \; .
\]
If $\oH^p_{\Fh} (X)$ is finite dimensional this agrees with the usual trace.
\end{defthm}

To motivate the definition of the trace in \ref{t32} further, assume that the cohomology spaces $\oH^p (\ep , H (\pi)^{\infty})$ are finite dimensional and choose bases in all of them. The union of these bases defines a basis of $\Hh$. In the situation 1) above writing $f^* : \Hh \to \Hh$ as a matrix in terms of these bases, we see that the sum of the diagonal entries is given by:
\[
\sum_{f^* \pi = \pi} \Tr (f^* \tei \oH^p (\ep , H (\pi)^{\infty})) \; .
\]
The more general case 2) can be similarly motivated.

\begin{remarknn}
 In example \ref{t36} below we will see that $\oH^p_{\Fh} (X)$ can be infinite dimensional even if $\Fh$ has dense leaves. 
\end{remarknn}

\begin{proofof}
  {\bf \ref{t32}} 1) Since $f : X \to X$ is a finite covering map we have
\[
f_* \mu = \mu \; .
\]
Note here that $f_* \mu$ is again $G$-invariant of volume one. Thus $f$ induces an isometric embedding
\[
f^* : L^2 (X) \longrightarrow L^2 (X) \; .
\]
In general this map is not $G$-equivariant. However because of the relation
\begin{equation}
  \label{eq:8}
  R_g \verk f^* = f^* \verk R_{f (g)} \quad \mbox{on} \; L^2 (X) \; \mbox{for all} \; g \in G \; ,
\end{equation}
it maps irreducible subrepresentations of $L^2 (X)$ into irreducible subrepresentations. Here isomorphic subrepresentations get mapped to isomorphic ones. Using (\ref{eq:8}) it follows that for every isotypical component $H (\pi)$ we have
\[
f^* H (\pi) = H (f^* \pi) \; .
\]
Equivalently:
\[
f^* (H (\pi)^{\infty}) = H (f^* \pi)^{\infty}
\]
under the pullback map $f^* : C^{\infty} (X) \to C^{\infty} (X)$. Thus the map $f^* = (f_* , f^*)$ on
\[
\oH^p_{\Fh} (X) = \oH^p_{\Fh} (\ep , C^{\infty} (X))
\]
permutes the subspaces $\oH^p (\ep , H (\pi)^{\infty})$ as stated. That $f^* \pi = \pi$ is equivalent to $\pi$ being the trivial representation is a special case of lemma \ref{t31}.

2) For endomorphisms $f : G \to G$ we do not have an induced map on $L^2 (X)$ unless $f$ is surjective. However the pullback map on $C^{\infty} (X)$ is all that is needed. Assume that $f^* (H (\pi)^{\infty})$ is not orthogonal to $H (\pi)^{\infty}$. Then
\[
P_{\pi} (f^* (H (\pi)^{\infty})) \neq 0 \; ,
\]
where $P_{\pi}$ denotes the orthogonal projection from $L^2 (X)$ onto $H (\pi)$. Since $P_{\pi}$ maps smooth vectors to smooth vectors, it follows that
\[
S = P_{\pi} \, |_{f^* (H (\pi)^{\infty})} : f^* (H (\pi)^{\infty}) \longrightarrow H (\pi)^{\infty}
\]
is a non-zero  $G$-intertwining map. So the map
\[
\tilde{S} = S \verk f^* : M_{\pi} \otimes V^{\infty}_{\pi} \cong H (\pi)^{\infty} \longrightarrow H (\pi)^{\infty} \cong M_{\pi} \otimes V^{\infty}_{\pi}
\]
is non-zero and satisfies
\[
\rho (g) \tilde{S} = \tilde{S} \rho (f (g))
\]
where $\rho = \rho^{\infty}_{\pi}$ is the representation on $V^{\infty}_{\pi}$. Note here that:
\[
R_g \verk f^* = f^* \verk R_{f (g)} \quad \mbox{on} \; C^{\infty} (X) \; \mbox{for all} \; G \; .
\]
Since $G$ acts trivially on the finite dimensional space $M_{\pi}$, it follows that there is a non-zero $G$-intertwining map $T$ from $V^{\infty}_{\pi}$ to itself satisfying
\[
\rho (g) T = T \rho (f (g)) \quad \mbox{for all} \; g \in G \; .
\]
In other words, $T$ is a non-zero $G$-intertwining map from $f^* (V^{\infty}_{\pi})$ to $V^{\infty}_{\pi}$. Applying lemma \ref{t31} it follows that $\pi$ is the trivial representation. The remaining assertions are clear.
\end{proofof}

\begin{theorem}
  \label{t34} For a generalized algebraic Anosov map $f$ for $(G , \Gamma)$ and a homogenous foliation $\Fh$ corresponding to $\ep$ as above, we have
\[
\sum^{\dim \Fh}_{p =0 } (-1)^p \Tr (f^* \tei \oH^p_{\Fh} (X)) = \sum_{x \atop f (x) = x} \varepsilon_x \, |\ddet (1- T_x f \tei T_x X / T_x \Fh)|^{-1} 
\]
where
\[
\varepsilon_x = \sgn \det (1 - T_x f \tei T_x \Fh) \; .
\]
\end{theorem}

\begin{proof}
  The proof of theorem \ref{t41_neu} applies verbally in this more general situation.
\end{proof}

Formula (\ref{eq:1}) of the introduction can now be extended. Let $f$ be a generalized algebraic Anosov diffeomorphism. Denote by $\eg_{\lambda}$ the generalized $\lambda$-eigenspace of $f_*$ on $\eg_{\C}$ and set
\[
\eg^u_{\C} = \bigoplus_{|\lambda|> 1} \eg_{\lambda} \; , \; \eg^s_{\C} = \bigoplus_{|\lambda|< 1} \eg_{\lambda} \; , \; \eg^e_{\C} = \bigoplus_{|\lambda|= 1} \eg_{\lambda} \; .
\]
Because of the relation
\[
[\eg_{\lambda} , \eg_{\mu}] \subset \eg_{\lambda\mu}
\]
these are sub Lie algebras of $\eg_{\C}$. For $\eg^* = \eg^*_{\C} \cap \eg$ we have $\eg^*_{\C} = \eg^* \otimes_{\R} \C$. Hence there is an $f_*$-invariant decomposition
\[
\eg = \eg^u \oplus \eg^s \oplus \eg^e \; .
\]
Here $\eg^s$ and $\eg^u$ can also be characterized as the set of $v \in \eg$ such that $f^N_* (v)$ tends to zero as $N \to + \infty$ resp. $N \to - \infty$. 

The subalgebras $\eg^s , \eg^u , \eg^e$ of $\eg$ are called the stable resp. unstable resp. neutral algebra of $f$. Let $\Fh^*$ be the foliation corresponding to $\eg^*$. Then we have
\[
TX = T \Fh^u \oplus T \Fh^s \oplus T\Fh^e
\]
and hence
\[
T_x X / T_x \Fh^u = T_x \Fh^s \oplus T_x \Fh^e \; .
\]
Using (\ref{eq:11}) it follows that for a fixed point $x$ of $f$ on $X$ we have
\[
T_x \Fh^s \oplus T_x \Fh^e = T^{se}_x \; , \; T_x \Fh^u = T^u_x
\]
where $T^*_x = T^*_{x\C} \cap T_x X$ and $T_{x\C}^{se}$ resp. $T^u_{x\C}$ is the sum of the eigenspaces of $T_xf$ on $T_x X \otimes \C$ corresponding to eigenvalues $\lambda$ of absolute value $|\lambda| \le 1$ resp. $|\lambda| > 1$. 

\begin{cor}
  \label{t35}
For every generalized algebraic Anosov diffeomorphism $f$ of $(G , \Gamma)$ the following formula holds:
\[
\sum_p (-1)^p \Tr (f^* \tei \oH^p_{\Fh^u} (X)) = \sum_{x \atop f (x) = x} \varepsilon_x \, |\ddet (1 - T_x f \tei T^{se}_x) |^{-1}
\]
where
\[
\varepsilon_x = \sgn \det (1 - T_x f \tei T^u_x) \; . 
\]
\end{cor}

The following example shows that even if the unstable foliation of a generalized algebraic Anosov diffeomorphism has dense leaves, its reduced leafwise cohomology can be infinite dimensional.

\begin{example}
  \label{t36}
\rm Let $G$ be the $3$-dimensional simply connected real Heisenberg group and let $\Gamma$ be its standard lattice. We can think of $G$ as the subgroup of matrices
\[
[x,y,z] := 
\begin{pmatrix}
  1 & x & z \\
0 & 1 & y \\
0 & 0 & 1
\end{pmatrix} \quad \mbox{in} \; \GL_3 (\R) \; .
\]
Then $\Gamma$ consists of those matrices $[n,m,k]$ with integer entries. The quotient $X = \Gamma \ohne G$ is compact. We will identify the Lie algebra $\eg$ of $G$ with the Lie algebra of matrices
\[
[|a,b,c|] := 
\begin{pmatrix}
  0 & a & c \\
0 & 0 & b \\
0 & 0 & 0 
\end{pmatrix} \quad \mbox{in} \; M_3 (\R) \; .
\]
One checks that $f : G \to G$ defined by
\[
f [x,y,z] = [x+y , 2x + y , -z + (x+y)^2 - \halb y (y-1)]
\]
is a Lie group automorphism with inverse:
\[
f^{-1} [x,y,z] = [y-x , 2x-y , x^2 - z - \halb (2x-y) (2x-y-1)] \; .
\]
The formulas show that $f (\Gamma) \subset \Gamma$ and $f^{-1} (\Gamma) \subset \Gamma$ hence $f (\Gamma) = \Gamma$. The differential $f_* : \eg \to \eg$ of $f$ is given by
\[
f_* [| a,b,c|] = [| a+b , 2a + b , \frac{b}{2} - c|] \; .
\]
Its eigenvalues are $1 + \sqrt{2} , 1 - \sqrt{2} , -1$. Hence $f$ is a generalized algebraic Anosov diffeomorphism. An eigenvector for $\lambda = 1 + \sqrt{2}$ is given by
\[
v = [| 1 , \sqrt{2} , 0 |] \; .
\]
Hence $\eg^u = \R \cdot v$. As $1 , \sqrt{2}$ are $\Q$-linear independent the one-dimensional foliation $\Fh^u$ has dense leaves \cite{AGH} Ch. IV and according to \cite{DS} Proposition 1.2 the reduced cohomology $\oH^1_{\Fh^u} (X)$ is {\it infinite dimensional}. Incidentally because of the relation
\[
\det (1 - f_* \tei \eg) = -4
\]
the automorphism $f$ of $X$ has $4$ different fixed points.
\end{example}

\begin{minipage}[t]{6cm}
A. Deitmar\\
School of mathematical sciences\\
University of Exeter\\
Laver building\\
North park road\\
Exeter EX4 4QE\\
Devon, UK\\
a.h.j.deitmar@exeter.ac.uk
\end{minipage} \hfill 
\begin{minipage}[t]{6cm}
C. Deninger \\
Mathematisches Institut\\
Einsteinstr. 62\\
48149 M\"unster\\
Germany\\
deninge@math.uni-muenster.de
\end{minipage}
\end{document}